\documentclass[draft, 12pt]{article}
\title{On the  Dimension of The   Virtually Cyclic Classifying
Space  of a Crystallographic Group  }
\author{Frank Connolly, Benjamin Fehrman and Michael Hartglass}
\date{\today}
\usepackage{amsmath,amsthm,amscd,amssymb}
\chardef\bslash=`\\ % p. 424, TeXbook
\theoremstyle{plain}

\theoremstyle{definition}

\theoremstyle{remark}

\newtheorem{thm}{Theorem}[section]

\newtheorem{Bieberbach1}{Theorem}[section]
\newtheorem{evcg}[Bieberbach1]{Definition}
\newtheorem{Espace}{Proposition}[section]
\newtheorem{Contractible}{Lemma}[section]
\newtheorem{UniqueC}[Contractible]{Claim}
\newtheorem{contrainfinite}[Contractible]{Claim}
\newtheorem{Liegroup}{Conjecture}[section]
\newcommand{\cF}{{\mathcal F}}
\newcommand{\cFin}{{\mathcal{FIN}}}

\newcommand{\cVC}{{\mathcal{VC} }}

\numberwithin{equation}{section}

\swapnumbers \theoremstyle{plain}

\begin{document}
%\date{\today}
%\subjclass{57 N 15,  57 R 67}
\maketitle
%%%%%%%%%%%%%%%%%%%%%%%%%%%%%%%%%%%%%%

\newcommand{\U}{\cal U \rm_{n}^{top}(X)}

%%%%%%%%%%%%%%

%%%%%%%%%%%%%%%%%%%%%%%
\begin{abstract}In this paper we construct a model for the classifying space, $B_{\mathcal{VC}}\Gamma$,
of a crystallographic group $\Gamma$ of rank $n$ relative to the
family $\mathcal{VC}$ of virtually-cyclic subgroups of $\Gamma$. The
model is used to show that there exists no other model for the
virtually-cyclic classifying space of $\Gamma$ with dimension less
than vcd($\Gamma$)+1, where vcd($\Gamma$) denotes the virtual
cohomological dimension of $\Gamma$.  In addition, the dimension of
our construction realizes this limit.\end{abstract}

\section{Introduction. Statement of Results}\

Let $\Gamma $ be a discrete  group.  Its classifying space, $B\Gamma
= E\Gamma/\Gamma$, is a polyhedron for which the $\Gamma$-space
$E\Gamma$ is a terminal object in the  category of    free
$\Gamma$-CW-complexes  and homotopy classes of $\Gamma$-maps.
$B\Gamma$ is therefore unique up to homotopy type.

In recent years topologists and others (e.g. Brown \cite{Brown},
Farrell and Jones \cite{FJ}, Lueck \cite{Lueck} and Serre
\cite{Serre}), have found great use for the more general notion of a
\emph{classifying space of a discrete group $\Gamma$,
 relative  to a family} $\cF$, of subgroups of $\Gamma$. This is a space
$B_{\cF}\Gamma= E_{\cF}\Gamma/ \Gamma$, where $E_{\cF}\Gamma$ is a
$\Gamma$-CW complex which is a terminal object in the category  of
homotopy classes of $\Gamma$-maps between $\Gamma$-CW-complexes
whose isotropy groups are in $\cF$.

According to Eilenberg and Ganea \cite{EG},  if $\Gamma$ has
cohomological dimension $\le n$, then one can construct a model of
$B\Gamma$ of dimension $\le max(n,3)$. There are similar results
showing $hdim (B_{\cFin}\Gamma)\le max(n,3)$, in many cases (see
\cite{CK}, \cite{Lueck}), but not all (see \cite{Brady}), whenever
$vcd(\Gamma)\le n$. Here $\cFin$ denotes the family of finite
subgroups of $\Gamma$, and
 we write
\[
hdim(Y) = min\{dim(X)| \text{ X is a CW-complex homotopy equivalent
to Y}\}
\]

A remarkable paper of F.T. Farrell and L. Jones \cite{FJ} introduces
the family
\[
\cVC = \{H|\; H \text{ is a virtually cyclic subgroup of } \Gamma\}
\]
(A group  is \emph{virtually cyclic} if it contains a cyclic
subgroup of finite index). They conjecture there that the $K$ or $L$
theory of a group $\Gamma$ can be computed from the homology of
$B_{\cVC} \Gamma $, taken with stratified coefficients in the $K$ or
$L$  theory of the virtually cyclic subgroups of $\Gamma$. The
K-theoretic version of this conjecture has been proved, for many
groups,  by Bartels and Reich \cite{Bart}.

Farrell and Jones \cite{FJ} also give two constructions of this
space $ B_{\cVC}\Gamma$, one of which  is a finite dimensional
CW-complex when $\Gamma$ is a discrete subgroup of a Lie Group.

The goal of this note is to provide a  geometrically simple
construction of $B_{\cVC}\Gamma$, as an $n+1$-dimensional
CW-complex, when $\Gamma$ is a crystallographic group of rank $n$
(see Section \ref{construct}). We use this construction to prove the
following result:

\

\begin{thm} Let $\Gamma$ be a crystallographic group of rank $n\geq2$. Then:
\[hdim(B_{\cVC}\Gamma) = n+1
\]
\end{thm}

Note that a crystallographic group $\Gamma$ of rank $1$ is
virtually-cyclic.  So in this case $hdim(B_{\mathcal{VC}}\Gamma)=0$.

\section{Basic Ideas.}\label{Basic}\

A crystallographic group is a discrete co-compact subgroup of
$\textrm{Iso}(\mathbb{R}^n)$, the group of isometries of
$\mathbb{R}^n$. In this paper, the \emph{translation subgroup} of a
crystallographic group $\Gamma$ will be
denoted\begin{displaymath}A=\Gamma\cap\textrm{Trans}(\mathbb{R}^n),\end{displaymath}
where $\textrm{Trans}(\mathbb{R}^n)$ denotes the subgroup of
$\textrm{Iso}(\mathbb{R}^n)$ consisting of translations.  $A$ is
normal in $\Gamma$. The \emph{point group} (or \emph{holonomy
group}) of $\Gamma$ will be written
\begin{displaymath}G=\Gamma/A.\end{displaymath}

The following theorem of Bieberbach [3] shows that $G$ is a finite
group.

\begin{Bieberbach1} Let $\Gamma$ be a crystallographic group
of rank $n$.  Then, $A$ is a finitely generated, free abelian group
of rank $n$ with finite index in $\Gamma$.\end{Bieberbach1}

A collection of subgroups $\mathcal{F}$ of a group $\Gamma$ is a
called a \emph{family} if $\mathcal{F}$ is closed under taking
subgroups and under conjugation in $\Gamma$.
\begin{evcg}\label{evcg} Let $\mathcal{F}$ be a family of subgroups
of $\Gamma$.  A a $\Gamma$-space $E$ is called
$\mathcal{F}$-universal if it satisfies the following
conditions:\begin{displaymath}\begin{array}{l} i.)\;\;\;E^H \textrm{
is contractible.}\;\forall
\;H\in\mathcal{F} \\
ii.)\;\;E^H=\emptyset.\;\forall \;H\notin\mathcal{F}
\\ iii.) \;E \textrm{ is a $\Gamma$-$CW$-complex}
\end{array} \end{displaymath}where $E^H$ denotes the fixed set of a subgroup
$H\subset\Gamma$ in $E$.\end{evcg}Then, we say $E/\Gamma$ is a
classifying-space for $\Gamma$ relative to $\mathcal{F}.$  By
\cite{Lueck} this specifies $E/\Gamma$ uniquely up to homotopy type.

We will mainly be interested in the
family:\begin{displaymath}\mathcal{VC}=\{H\subset\Gamma\;|\;H\textrm{
a virtually-cyclic subgroup}\}.\end{displaymath}  Later, our model
for a $\mathcal{VC}$-universal space of a crystallographic group
$\Gamma$ will be denoted $E_\mathcal{VC}\Gamma$.  The classifying
space for $\Gamma$ relative to $\mathcal{VC}$ will be written
$B_\mathcal{VC}\Gamma=E_\mathcal{VC}\Gamma/\Gamma$.

Finally, let $\mathcal{C}=\{H\subset A\;|\;H\textrm{ is a maximal
cyclic subgroup of A}\}$ be the set of maximal cyclic subgroups of
$A$. For each subgroup $C\in\mathcal{C}$ we
define:\begin{displaymath}
\mathbb{R}^{n-1}(C)=\{l\subset\mathbb{R}^n\;|\;l\textrm{ a line,
}C\cdot l=l\}.\end{displaymath} This is the set of lines in
$\mathbb{R}^n$ left invariant under the group action of $C$.  The
quotient map $\pi_C:\mathbb{R}^n\rightarrow\mathbb{R}^{n-1}(C)$ is
the map:\begin{displaymath}\pi_C(x)=\textrm{ the unique line
}l\in\mathbb{R}^{n-1}(C)\textrm{ containing } x\;\;\forall\;
x\in\mathbb{R}^n.\end{displaymath}

The Hausdorff metric on the non-empty closed sets of $\mathbb{R}^n$
restricts to a metric on $\mathbb{R}^{n-1}(C)$. In addition, the
quotient topology on $\mathbb{R}^{n-1}(C)$ coincides with the metric
topology on $\mathbb{R}^{n-1}(C)$.  Note that $\mathbb{R}^{n-1}(C)$
is isometric to $\mathbb{R}^{n-1}$ in such a way that $\pi_C$ is a
linear map.

\section{Construction of $E_{\mathcal{VC}}\Gamma$}\label{construct}\

Let $\Gamma$ be a crystallographic group of rank $n$, with holonomy
group $G$. Our model for $B_{\mathcal{VC}}\Gamma$ will be the orbit
space of a $G$-action on an infinite union of solid $n+1$-tori
sharing a common boundary. $E_{\mathcal{VC}}\Gamma$ will be a
similar union of mapping cylinders sharing a common source, one
cylinder for each $C\in\mathcal{C}$.

Equip $\mathcal{C}$ with the discrete topology.  Define an
equivalence relation $\sim$ on $\mathbb{R}^n\times
I\times\mathcal{C}$ as follows: $(x,t,C)\sim(x',t',C')$ if
\begin{displaymath}\begin{array}{l} i.)\;\;\;0<t=t'<1,\;C=C'\textrm{
and }x=x',\textrm{ or} \\ ii.)\;\;t=t'=1,\;C=C' \textrm{ and
}\pi_C(x)=\pi_C(x'),\textrm{ or}
 \\ iii.)\;t=t'=0\textrm{ and }x=x'.\end{array}\end{displaymath}
We define:
\begin{equation}\label{E}E=(\mathbb{R}^n\times
I\times\mathcal{C})/\sim\end{equation} with the quotient topology.
The equivalence class of an element $(x,t,C)$ is written $[x,t,C]$.
For each $C\in\mathcal{C}$ we
write\begin{displaymath}Cyl(\pi_c)=\left(\mathbb{R}^n\times
I\times\{C\}\right)/\sim.\end{displaymath}Note this subspace of $E$
is just the mapping cylinder of
$\pi_C:\mathbb{R}^n\rightarrow\mathbb{R}^{n-1}(C)$, and $E$ is the
union of these subspaces $Cyl(\pi_C)$.

The action of $\Gamma$ on $E$ is defined
by:\begin{equation}\label{action} \gamma\cdot[x,t,C]=[\gamma \cdot
x,t,\gamma C\gamma^{-1}]\;\forall\; \gamma\in\Gamma.\end{equation}

In the following sections we prove:

\begin{Espace}\label{Espace} The $\Gamma$-space $E$ constructed in (\ref{E}) and
(\ref{action}) is a $\mathcal{VC}$-universal $\Gamma$-space (see
Definition \ref{evcg}).\end{Espace}

\section{$E$ is a $\Gamma$-CW-Complex}\label{CW}\

In this section, we show that $E$ has the structure of a
$\Gamma$-CW-complex.  Choose any $\Gamma$-CW-complex structure on
$\mathbb{R}^{n}$, called $X$, whose cells, $e$ $\in$ $Cell(X)$ are
all convex polytopes.  Here, $Cell(X)$ is the collection of cells of
$X$. It is enough to show that for each $C$ $\in$ $\mathcal{C}$,
this CW-structure extends equivariantly to a CW-structure on
$Cyl(\pi_{C})$.  This would be obvious if $\pi_{C}$: $X$
$\longrightarrow$ $\mathbb{R}^{n-1}(C)$ were cellular with respect
to a CW-complex structure $Y$ on $\mathbb{R}^{n-1}(C)$, for we would
then use the resulting CW-complex structure on $Cyl(\pi_{C})$.  It
is nearly as obvious if there is a subdivision $X'$ of $X$ (that is,
each cell of $X'$ will be a subset of a cell of $X$) relative to
which $\pi_{C}$: $X'$ $\longrightarrow$ $Y$ is cellular, since $X
\times I$ can be subdivided by replacing only the cells of $X \times
\{1\}$ with the cells of $X' \times \{1\}$ and leaving the other
cells unchanged.  This means that it suffices to show that
$\mathbb{R}^{n-1}(C)$ has the structure of an
$N_\Gamma(C)$-CW-complex so that $\pi_{C}$ is cellular relative to
the subdivision $X'$ of $X$.

First note \{$\pi_{C}(e)$ $\mid$ $e$ $\in$ $Cell(X)$\} is a locally
finite collection in $\mathbb{R}^{n-1}(C)$, because $Cell(X)$ is a
locally finite collection in $\mathbb{R}^{n}$.  We note that a
polytope generated by $n$ points is the convex hull of those $n$
points, and we choose $X$ so that $Cell(X)$ is generated by convex
polytopes. For each $y$ $\in$ $\mathbb{R}^{n-1}(C)$ define:
\begin{displaymath}f_{y} = \bigcap \{\pi_{C}(e) \mid y \textrm{ is
contained in the interior of the convex polytope }
\pi_{C}(e)\}.\end{displaymath}  This intersection is finite, so we
have \begin{displaymath}Int(f_{y}) = \bigcap \{ Int(\pi_{C}(e)) \mid
y \in Int(\pi_{C}(e))\}.\end{displaymath} Note $y$ $\in$
$Int(f_{y})$ so each point of $\mathbb{R}^{n-1}(C)$ is in the
interior of exactly one of the convex polytopes $f_{y}$. Also,
\{$f_{y}$ $\mid$ $y$ $\in$ $\mathbb{R}^{n-1}(C)$\} is a locally
finite collection.  To see that these form an
$N_\Gamma(C)$-CW-complex we show that if $dim(f_{y})$ = $d$, then
$\dot f_{y}$ is a union of cells of dimension $<$ $d$ (here, $\dot
f_{y}$ denotes the boundary of $f_{y}$).  Let $y'$ $\in$ $\dot
f_{y}$.  Then $y'$ $\in$ $\pi_{C}(e)^{\cdot}$ for some $e$ such that
$y$ $\in$ $Int(\pi_{C}(e))$.  Note that $ \pi_{C}(e)^{\cdot}$ is a
subcomplex of $\pi_{C}(\dot e)$. But $dim(\pi_{C}(e))$ = $d$, so
$\pi_{C}(e)^{\cdot}$ is a union of cells $\pi_{C}(e')$ of dimension
$\leq$ $d-1$.  Hence, there is a cell $\pi_{C}(e')$ with $y'$ $\in$
$Int(\pi_{C}(e'))$ such that $\pi_{C}(e')$ $\subset$
$\pi_{C}(e)^{\cdot}$. Therefore, $\dot f_{y}$ is a union of cells of
dimension $<$ $d$. This shows that $Cell(Y)$ := \{$f_{y}$ $\mid$ $y$
$\in$ $\mathbb{R}^{n-1}(C)$\} is the collection of cells of an
$N_\Gamma(C)$-CW-complex structure on $\mathbb{R}^{n-1}(C)$ denoted
$Y$.  The cells of $X'$ are now easily defined.  They are the
nonempty sets of the form $e$ $\bigcap$ $\pi_{C}^{-1}(f)$ with e
$\in$ $Cell(X)$ and $f$ $\in$ $Cell(Y)$. Note that $e$ $\bigcap$
$\pi_{C}^{-1}(f)$ is a convex polytope contained in the cell $e$.
Since $\pi_{C}$ takes each cell into a cell linearly, it follows
that $\pi_{C}$: $X'$ $\longrightarrow$ $Y$ is cellular. The above
observations show that $E$ admits the structure of a
$\Gamma$-CW-complex.

\section{Contractibility of Fixed Sets}\

In this section we prove:

\begin{Contractible}\label{Contractible} Let $E$ be the $\Gamma$-space
constructed in (\ref{E}).  For each subgroup $H\subset\Gamma$, $E^H$
is contractible if $H\in\mathcal{VC}$ and $E^H$ is empty if
$H\notin\mathcal{VC}$.

\begin{proof}  Let $H$ be a subgroup of $\Gamma$. Then,
\begin{equation}\label{finitefix}(\mathbb{R}^n)^H\textrm{ is }\left\{
\begin{array}{ll} \textrm{contractible} & \textrm{if } |H|<\infty \\ \emptyset
& \textrm{if }|H|=\infty.\end{array} \right. \end{equation}

The $\Gamma$-action of (\ref{action}) defines an action of
$N_\Gamma(C)/C$ on $\mathbb{R}^{n-1}(C)$.  This gives a map
$j_C:N_\Gamma(C)/C\rightarrow\textrm{Iso}(\mathbb{R}^{n-1}(C))$.
$\textrm{Ker}(j_C)$ has at most two elements. Recall
$\mathbb{R}^{n-1}(C)$ is isometric to $\mathbb{R}^{n-1}.$  The image
of $j_C$ is a crystallographic group.

Assume $H$ is a subgroup of $N_\Gamma(C)$ for some
$C\in\mathcal{C}$. By the above,
\begin{equation}\label{infcontra} (\mathbb{R}^{n-1}(C))^H\textrm{ is }\left\{
\begin{array}{ll} \textrm{contractible} & \textrm{if }|HC/C|<\infty \\
\emptyset & \textrm{if }|HC/C|=\infty.\end{array} \right.
\end{equation}

\begin{UniqueC}Let
$H\subset\Gamma$ satisfy $|H|=\infty$.  Then, there exists at most
one $C\in\mathcal{C}$ satisfying $E^H\cap Cyl(\pi_C)\neq\emptyset$.

\begin{proof}Suppose there are two subgroups $C,C'\in\mathcal{C}$
satisfying $E^H\cap Cyl(\pi_C)\neq\emptyset$ and $E^H\cap
Cyl(\pi_{C'})\neq\emptyset$.  By (\ref{action}) this implies
$H\subset N_\Gamma(C)$ and $H\subset N_\Gamma(C')$. Then, by
(\ref{finitefix}) and (\ref{infcontra}) we know $|HC/C|<\infty$ and
$|HC'/C'|<\infty$. Therefore, the relation $HC/C\approx H/(H\cap C)$
shows that $H\cap C$ and $H\cap C'$ have finite index in $H$. So
$H\cap C\cap C'$ has finite index in $H$.  But $H\cap C\cap
C'=\{e\}$ and $|H|=\infty$, a contradiction.
\end{proof}\end{UniqueC}

\begin{contrainfinite}\label{contrainfinite}Let $H\subset\Gamma$ be an infinite
subgroup.  Then $E^H$ is contractible if $H\in\mathcal{VC}$ and
$E^H$ is empty if $H\notin\mathcal{VC}$.

\begin{proof}If $H\in\mathcal{VC}$ there exists precisely one
$C\in\mathcal{C}$ satisfying $E^H\cap Cyl(\pi_C)\neq\emptyset$
(namely, the unique element of $\mathcal{C}$ containing $H\cap A$).
For this $C$, (\ref{finitefix}) implies
$E^H=(\mathbb{R}^{n-1}(C))^H\times\{1\}\times\{C\}$. This space is
contractible by (\ref{infcontra}). Conversely, suppose there exists
a unique $C\in\mathcal{C}$ satisfying $E^H\cap
Cyl(\pi_C)\neq\emptyset$. Then, by (\ref{finitefix}) and
(\ref{infcontra}) we know $H\in\mathcal{VC}$ because $|H/(H\cap
C)|<\infty$. Therefore, $E^H$ is contractible. Finally, if there
exists no $C\in\mathcal{C}$ satisfying $E^H\cap
 Cyl(\pi_C)\neq\emptyset$ then
$E^H=\emptyset$.\end{proof}\end{contrainfinite}

Thus, for a subgroup $H\subset\Gamma$ the set $E^H$ is the
following:

\small

\begin{equation}\label{Econtra}E^H=\left\{\begin{array}{ll}
(\mathbb{R}^n)^H\cup (\cup\{Cyl(\pi_C^H)\;|\;H\subset N_\Gamma(C)\})
&
|H|<\infty \\
(\mathbb{R}^{n-1}(C))^H & |H|=\infty,\;H\in\mathcal{VC},\;H\subset
N_\Gamma(C) \\ \emptyset &
|H|=\infty,\;H\notin\mathcal{VC}.\end{array}\right.\end{equation}

\normalsize

If $H\in\mathcal{VC}$ and $H$ is infinite then $E^{H}$ is
contractible by Claim \ref{contrainfinite}.  If $H$ is finite, then
$\pi_C^{H}$ is a homotopy equivalence for each $C\in\mathcal{C}$ for
which $H\subset N_\Gamma(C)$ (here, $\pi_C^H$ is the restriction of
$\pi_C$ to $(\mathbb{R}^n)^H$). Therefore the subcomplex
$(\mathbb{R}^{n})^{H}$ is a strong deformation retract of the
CW-complex $Cyl(\pi_C^H)$, and these deformations coalesce to make
$(\mathbb{R}^{n})^{H}$ a strong deformation retract of $E^H$. This
and (\ref{finitefix}) show that $E^H$ is
contractible.\end{proof}\end{Contractible}

\begin{proof}(\emph{of Proposition \ref{Espace}}):  This is clear from Section
\ref{CW} and Lemma \ref{Contractible}.\end{proof}

Henceforth we denote the $\mathcal{VC}$-universal $\Gamma$-space $E$
by:
\begin{displaymath}E=E_{\mathcal{VC}}\Gamma\end{displaymath}and
\begin{displaymath}E/\Gamma=B_{\mathcal{VC}}\Gamma.\end{displaymath}

\section{Geometric Properties of $B_{\mathcal{VC}}\Gamma$}\

Let $\Gamma$, $A$ and $G$ be as in Section \ref{Basic}. In this
section we will show that $B_{\mathcal{VC}}\Gamma$ is the orbit
space of a $G$-action on an infinite union of solid $n+1$-tori
sharing a common boundary.

First we exhibit a
homeomorphism\begin{equation}\label{homeof}f:Cyl(\pi_C)/A\rightarrow
D^2\times T^{n-1}.\end{equation}  This will show $Cyl(\pi_C)/A$ is
homeomorphic to the solid $n+1$-torus $D^2\times T^{n-1}$ for each
$C\in\mathcal{C}$. Choose a basis \{$a_{1}, .., a_{n}$\} for $A$ so
that $a_{1}$ generates $C$. Set $v_{i}$ = $a_{i}(0)$. Then \{$v_{1},
..., v_{n}$\} is a basis for $\mathbb{R}^{n}$.  If $p$ = $x_{1}v_{1}
+ ... + x_{n}v_{n}$, we will write $p$ = $<x_{1}, ..., x_{n}>$.

Let $f:Cyl(\pi_C)/A\rightarrow D^2\times T^{n-1}$ be the
map:\begin{displaymath}f([[p,t,C]])=(1-t)e^{2\pi ix_1}\times e^{2\pi
ix_2}\times\dots\times e^{2\pi ix_n},\end{displaymath}where
$p=<x_{1},...,x_{n}>$ and $[[p,t,C]]$ denotes the orbit, in
$Cyl(\pi_C)/A$, of a point $[p,t,C]\in Cyl(\pi_C)$. The map $f$ is
continuous and bijective. It is a homeomorphism as $Cyl(\pi_C)/A$ is
compact.

Because
$E_\mathcal{VC}\Gamma/A=\cup\{Cyl(\pi_c)/A\;|\;C\in\mathcal{C}\}$
the space $E_\mathcal{VC}\Gamma/A$ is homeomorphic to an infinite
union of solid $n+1$-tori sharing the common boundary
$\mathbb{R}^n/A$. Note that $\mathbb{R}^n/A$ is homeomorphic to
$T^{n}$.

The relation $E_{\mathcal{VC}}\Gamma/\Gamma\simeq
(E_{\mathcal{VC}}\Gamma/A)/(\Gamma/A)$ proves\begin{displaymath}
B_{\mathcal{VC}}\Gamma=(E_{\mathcal{VC}}\Gamma/A)/G.\end{displaymath}
$B_{\mathcal{VC}}\Gamma$ is therefore the orbit space of a
$G$-action on a union of solid $n+1$-tori sharing a common boundary,
one torus for each $C\in\mathcal{C}$.

\section{Computation of $hdim(B_{\mathcal{VC}}\Gamma)$}\

In this section we prove Theorem 1.1.\\

\begin{proof} Let $i_{C}$: $Cyl(\pi_{C})/A$ $\rightarrow$
$E_{\mathcal{VC}}\Gamma/A$ be inclusion, $\tau$:
$E_{\mathcal{VC}}\Gamma/A$ $\rightarrow$ $B_{\mathcal{VC}}\Gamma$ be
the obvious projection and $\iota_{C}$ = $\tau$ $\circ$ $i_{C}$:
$Cyl(\pi_{C})/A$ $\rightarrow$ $B_{\mathcal{VC}}\Gamma$. We will
make a computation in the $n+1st$ homology group
$H_{n+1}(B_{\mathcal{VC}}\Gamma)$ using the cellular chain complex,
$C_{*}(B_{\mathcal{VC}}\Gamma)$. By (\ref{homeof}), for each $C$
$\in$ $\mathcal{C}$ there exists an n+1 chain, $\chi_{C}$ $\in$
$C_{n+1}(Cyl(\pi_{C})/A)$ such that $\partial_{n+1}(\chi_{C})$ =
[$T_{n}$] (where $\partial_{n+1}$ is the boundary map and [$T_{n}$]
is the fundamental cycle of $\mathbb{R}^{n}/A$). Choose $C$, $C'$
$\in$ $\mathcal{C}$ such that $C$ and $C'$ are not conjugate in
$\Gamma$ (here is where we use that $rank(\Gamma) \geq 2$).  Denote
$\iota_{C\#}$ and $\iota_{C'\#}$ as the chain maps induced from
$\iota_{C}$ and $\iota_{C'}$ respectively. Then, in
$C_{n+1}(B_{\mathcal{VC}}\Gamma)$:
\begin{displaymath}z= \iota_{C\#}(\chi_{C}) -
\iota_{C'\#}(\chi_{C'}) \neq 0. \textrm{ and }
\partial_{n+1}(z)=0.\end{displaymath} Indeed, \begin{eqnarray*}\partial_{n+1}(z) & = &
\partial_{n+1}\iota_{C\#}(\chi_{C}) -
\partial_{n+1}\iota_{C'\#}(\chi_{C'})\nonumber \\ & = &
\iota_{C\#}\partial_{n+1}(\chi_{C}) -
\iota_{C'\#}\partial_{n+1}(\chi_{C'}) = \iota_{C\#}[T_{n}] -
\iota_{C'\#}[T_{n}] = 0\nonumber\end{eqnarray*} as the restrictions
of $\iota_{C}$ and $\iota_{C'}$ to $\mathbb{R}^{n}/A$ are the same.
Since $C_{n+2}(B_{\mathcal{VC}}\Gamma)$ is trivial, we conclude that
$H_{n+1}(B_{\mathcal{VC}}\Gamma$) is nontrivial, implying that
$hdim(B_{\mathcal{VC}}\Gamma$) $\geq n+1$. But
$B_{\mathcal{VC}}\Gamma$ has dimension $n+1$. Therefore
$hdim(B_{\mathcal{VC}}\Gamma$) = $n+1$.\end{proof}

\section{Conclusion and an Open Question}

In this paper, we were able to construct a model for both
$E_{\mathcal{VC}}\Gamma$ and $B_{\mathcal{VC}}\Gamma$ for a rank-$n$
crystallographic group $\Gamma$.  Despite the fairly straightforward
constructions of our models, they are not locally finite
CW-complexes.  Therefore they can not be imbedded in any Euclidian
space as they are not metrizable.

One question that is still unsolved is the
following:\begin{Liegroup}Suppose $\Gamma$ is a discrete subgroup of
a Lie Group and $\Gamma$ is not virtually-cyclic.  Then,
$hdim(B_{\mathcal{VC}}\Gamma)=vcd(\Gamma)+1$.\end{Liegroup} We have
shown this is so if $\Gamma$ is a crystallographic group, of rank
$n\geq2$.

\pagebreak{}

\end{document}